\documentclass{article}
\usepackage{amsmath}
\usepackage{amssymb}
\usepackage{latexsym}
\usepackage{amsthm}
\usepackage{mathrsfs}
\usepackage{wasysym}
\usepackage{fancyhdr}
\usepackage{amsfonts}
\usepackage{amssymb}
\usepackage{epsfig}

\begin{document}
\title{Dominating Plane Triangulations} 
\author{Michael D. Plummer\thanks{
Department of Mathematics,
Vanderbilt University,
Nashville, TN 37215,
Email: \tt{michael.d.plummer@vanderbilt.edu}}\,,
Dong Ye\thanks{Department of Mathematical Sciences,
Middle Tennessee State University,
Murfreesboro, TN 37132, USA,
Email: \tt{dong.ye@mtsu.edu}}\;  and
Xiaoya Zha\thanks{Department of Mathematical Sciences,
Middle Tennessee State University,
Murfreesboro, TN 37132, USA, Email: {\tt{xiaoya.zha@mtsu.edu}}. Research supported by NSA Grant H98230-1-02192.}}

\date{}
\maketitle
\begin{abstract}   
In 1996, Tarjan and Matheson proved that if $G$ is a plane triangulated disc with $n$ vertices, $\gamma (G)\le n/3$,
where $\gamma (G)$ denotes the domination number of $G$.
Furthermore, they conjectured that the constant $1/3$ could be improved to $1/4$
for sufficiently large $n$.
Their conjecture remains unsettled.\par
In the present paper, it is proved that if $G$ is a hamiltonian plane triangulation with $|V(G)|=n 
$ vertices and
minimum degree at least 4, then $\gamma (G)\le\max\{\lceil 2n/7\rceil, \lfloor 5n/16\rfloor\}$.
It follows immediately that if $G$ is a 4-connected plane triangulation with $n 
$ vertices, then
$\gamma (G)\le\max\{\lceil 2n/7\rceil, \lfloor 5n/16\rfloor\} $.
It then follows that if $n\ge 26$,
then $\gamma (G)\le \lfloor 5n/16\rfloor$.\medskip

\end{abstract}

\section{  Introduction and Terminology}  
In 1996, Matheson and Tarjan [MT] proved that if $G$ is a plane triangulated disc then
$\gamma (G)\le |V(G)|/3$.
(In particular, then, the same bound on $\gamma$ applies to any triangulation of the plane.)
Plummer and Zha [PZ] proved that if $G$ is a triangulation of the projective plane, then $\gamma (G)\le |V(G)|/3$
and if $G$ is a triangulation of either the torus or Klein bottle, then $\gamma (G)\le \lceil |V(G)|/3\rceil$.
The latter result was sharpened by Honjo et al. [HKN] who showed that $\gamma (G)\le |V(G)|/3$
for graphs embedded in these two surfaces.
They also showed that for any surface $\Sigma$, there is a positive integer $\rho (\Sigma )$ such that
if $G$ is embedded as a triangulation of $\Sigma$ and the embedding has face-width at least $\rho (\Sigma )$,
then $\gamma (G)\le |V(G)|/3$.
Even more recently, Furuya and Matsumoto [FM] generalized this result by showing that $\gamma (G)\le |V(G)|/3$,
for every triangulation $G$ of any closed surface.\par

In [MT], Matheson and Tarjan also conjectured that their bound of $n/3$ could be improved; namely, if $n$ is sufficiently large, then any triangulation of the plane with $n$ vertices would have $\gamma\le n/4$.
The 
triangle has $\gamma (K_3)=1 = n/3$, the octahedron shown in Figure 1.1(i) has $\gamma=2=n/3$ and the 7-vertex graph shown in Figure 1.1(ii) has $\gamma = 2/7>1/4$, so this shows that
one must assume $n\ge 8$ in order for the Matheson-Tarjan conjecture to be true.

\begin{figure}[!hbtp] 
\begin{center}
\includegraphics[scale=.75]{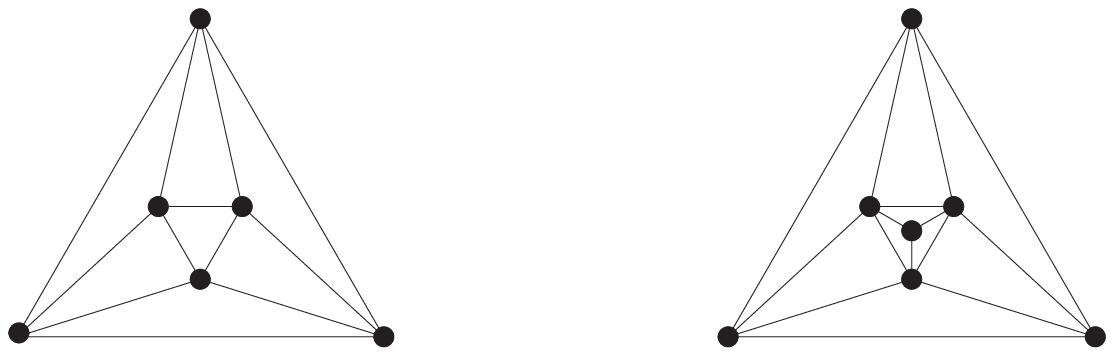} 
\centerline{\bf Figure 1.1}
\end{center}
\end{figure}

More generally, in [PZ] the first and third authors of the present paper conjectured that if $G$ is a
triangulation of {\it any}
non-spherical surface, then
$\gamma (G)\le n/4$.
Both these conjectures involving the $n/4$ bound remain unsettled.\par


In 2010, King and Pelsmajer [KP] proved the Matheson-Tarjan bound of n/4 holds in the plane 
case when the maximum degree of the triangulation is 6.\par

An outerplanar graph is a graph embedded in the plane in such a way that all vertices of the graph lie on the boundary of the external face.
An outerplanar graph is {\it maximal} (outerplanar) if it is not possible to add any new edge to $G$ without destroying outerplanarity.
In 2013, Campos and Wakabayashi [CW] proved that if $G$ is a maximal outerplanar graph with at least $n\ge 4$ vertices, then $\gamma (G)\le (n+t)/4$, where $t$ is the number of vertices of degree 2 in $G$.


It was proved independently by N\"unning [N] and by Sohn and Yuan [SY] that for any graph $G$ with
$n$ vertices and minimum degree $\delta (G)\ge 4$, $\gamma (G)\le 4n/11$.
In the present paper we will show that if $G$ is a plane triangulation on $n
$ vertices which has a Hamilton cycle and $\delta (G)\ge 4$, then $\gamma (G)\le \max\{\lceil 2n/7\rceil, \lfloor 5n/16\rfloor\}$.
It follows immediately that if $G$ is a 4-connected plane triangulation on $n
$ vertices, then
 $\gamma(G)\le \max\{\lceil 2n/7\rceil, \lfloor 5n/16\rfloor\}$.
Note that if $n\notin \{6,8,9,11,12,15,19,22,25\}$, $\lceil 2n/7\rceil \le \lfloor 5n/16\rfloor$, so this result can be restated to
say that if $G$ is a 4-connected plane triangulation on $n$ vertices and $n\ge 26$, then $\gamma (G)\le \lfloor 5n/16\rfloor$.

\section{Preferred Hamilton cycles in plane triangulations}

Let $G$ be a plane triangulation with $\delta (G)\ge 4$ and suppose $G$ contains a Hamilton cycle $H$.
We can think of $H$ bounding a triangulated inner subgraph $G_{int}$ and a triangulated outer subgraph
$G_{ext}$ such that $G_{int}\cap G_{ext}= H$.
Suppose $v\in V(G)$.
We denote by $i\deg v$ (respectively, $o\deg v$) the degree of vertex $v$ in $G_{int}$ (resp. in $G_{ext}$).
\par

We will need to pay particular attention to those vertices which have $i\deg v = 2$ or $o\deg v=2$.
We will call these vertices {\it 2-vertices}.
Examples are shown in Figure 2.1 where $i\deg v_1=2$ and $o\deg v_2=2$ respectively and hence each is an example of
a 2-vertex.\par

\begin{figure}[!hbtp]\refstepcounter{figure}
\begin{center}
\includegraphics[scale=.45]{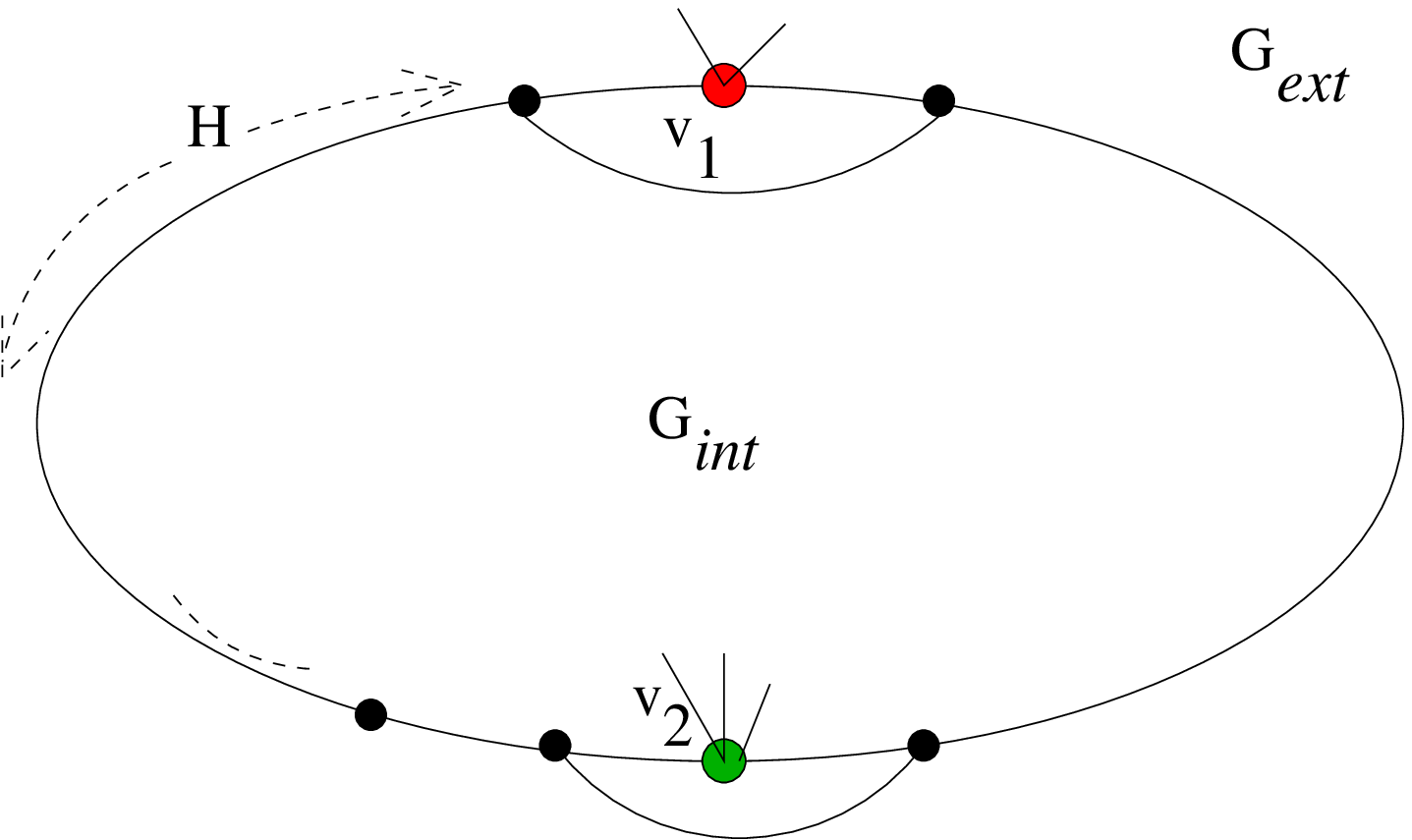} 
\centerline{\bf Figure 2.1}
\end{center}
\end{figure}

Let $G$ be a Hamiltonian plane triangulation.  Such a graph may have many Hamilton cycles.
We now show how to select such a cycle with certain properties we shall need later on.\par

To this end, let begin with any Hamilton cycle $H$ in $G$.
A triangle $T$ of $G_{int}$ will be called {\it internal} (with respect to $H$) if $E(T)\cap E(H)=\emptyset$.
\par
\medskip
\noindent
{\bf Lemma 2.1:}  Let $G$ be a plane triangulation with $\delta (G)\ge 4$ which contains a Hamilton cycle.
Suppose $\gamma (G)\ge 2$.
Then there exists a Hamilton cycle in $G$ containing no three consecutive 2-vertices.\par
\medskip\noindent
{\bf Proof:}  First, suppose that $G$ contains four consecutive 2-vertices.
Since $G$ is a triangulation, these four 2-vertices must alternate between $ideg=2$ and $o\deg=2$.
So let us suppose that we have four consecutive 2-vertices forming a subpath of $H$, which we will denote by
$wxyz$, such that 
$i\deg w = i\deg y = 2$ and $o\deg x = o\deg z=2$.
Let the predecessor of $w$ be $a$ and the successor of $z$ be $b$.
Note that $a\ne b$, for if $a=b$, then $G=K_5$ contradicting the hypothesis that $G$ is planar.

\begin{figure}[!hbtp] 
\begin{center}
\includegraphics[scale=.4]{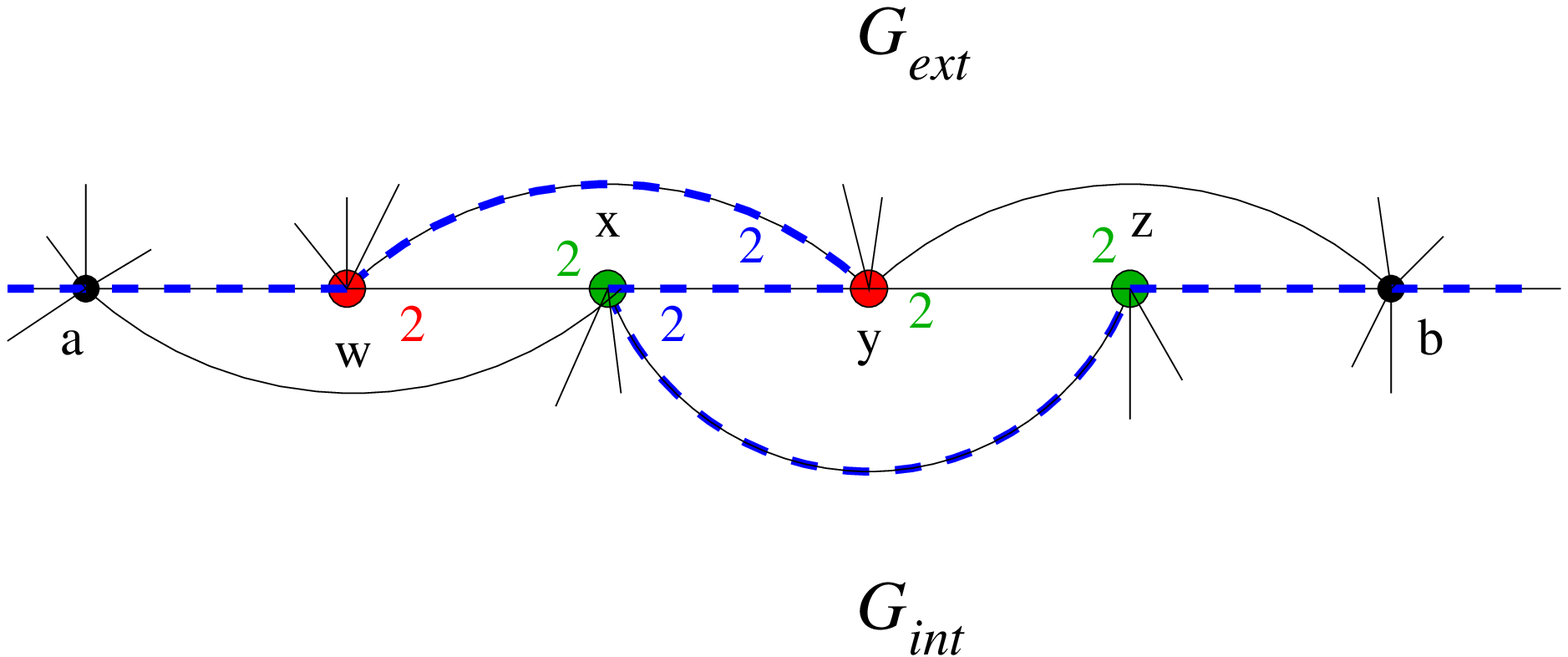} 
\centerline{\bf Figure 2.2}
\end{center}
\end{figure}

Form a new Hamilton cycle $H'$ substituting for the path $awxyzb$ the path $awyxzb$ as indicated by the dashed path shown in Figure 2.2.
With respect to the original Hamilton cycle $H$, $w,x,y$ and $z$ are 2-vertices.
However, with respect to $H'$, only $x,y$ and (possibly) $b$ are 2-vertices.
(Note that the vertices $w$ and $z$ are {\it not} 2-vertices with respect to $H'$ because the degrees of $w$ and$z$ are both at least 4.)\par

Next suppose that $H$ does not contain four consecutive 2-vertices, but suppose it does contain three consecutive 2-vertices.
Let three such consecutive 2-vertices be, in order, $x,y$ and $z$.
Let the predecessor of $x$ be $a$ and the successor of $z$ be $b$.
Again, $a\ne b$ as before.
Let the neighbor of $b$ on $H$ different from $z$ be denoted as $b'$.\par

\begin{figure}[!hbtp] 
\begin{center}
\includegraphics[scale=.35]{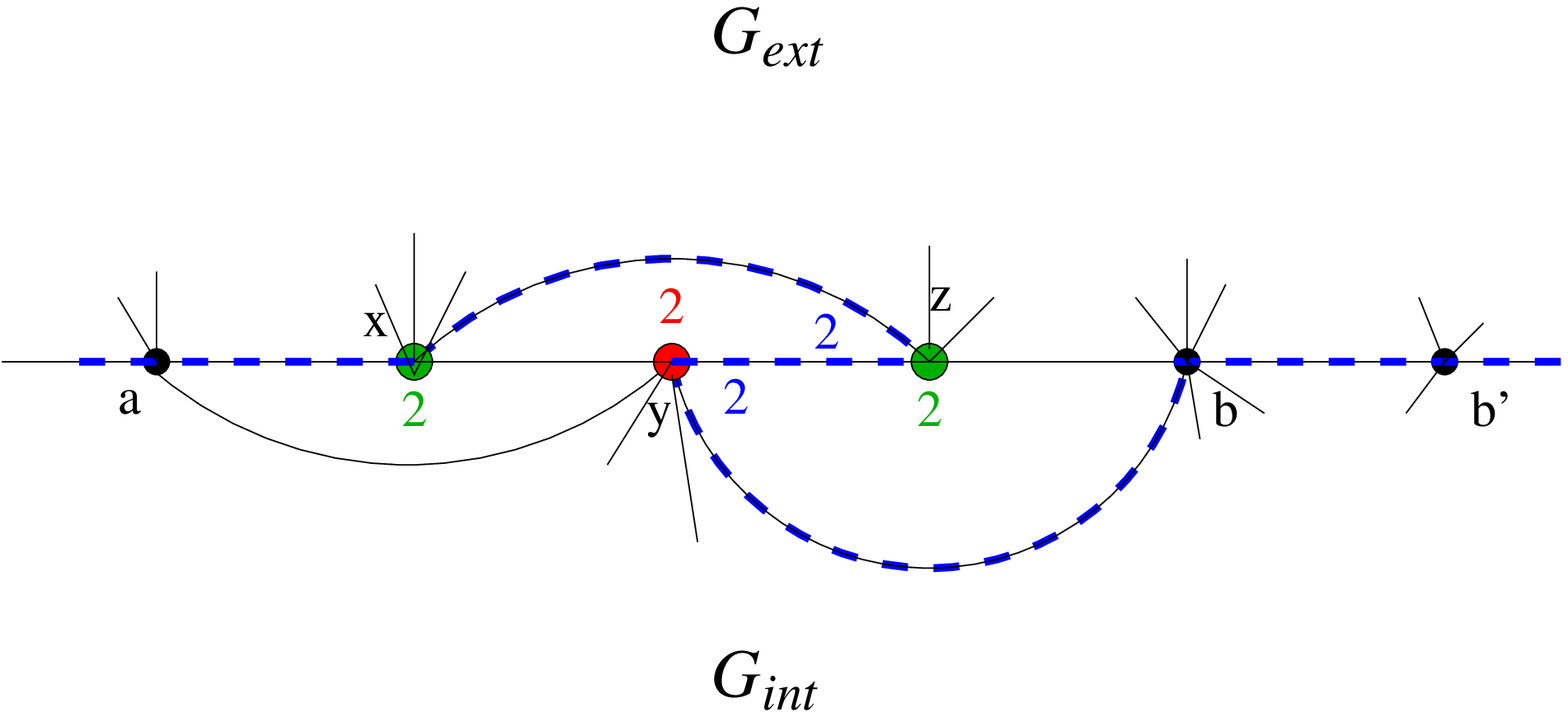} 
\centerline{\bf Figure 2.3}
\end{center}
\end{figure}

Now suppose $b$ has another neighbor via an edge in $G_{int}$ other than $y,z$ and $b'$.
Then replace the path $axyzb$ by $axzyb$ in $H$ to obtain a new Hamilton cycle $H'$.
Then whereas $H$ has $x,y$ and $z$ as 2-vertices, $H'$ has vertices $y$ and $z$
as 2-vertices.
But neither $x$ and $b$ is a 2-vertex with respect to $H'$ since the degree of $x$ is at least 4 and because of our assumption on the fourth neighbor of $b$.

Thus $H'$ has one less triple of successive 2-vertices than does $H$.
(See Figure 2.3.)

Suppose, on the other hand, that the only neighbors of vertex $b$ via edges in $G_{int}$ are $y,z$ and $b'$.
Hence $y$ is adjacent to $b'$.
In this case, vertex $b$ will be a new 2-vertex with respect to $H'$, where $H'$ is obtained from $H$ be replacing the path $axyzb$
by the path $axzyb$ and hence the total number of 2-vertices in $H'$ is the same as in $H$.

More generally, suppose that $y$ is adjacent to $b'=b_0, b_1,\ldots,b_k,b_{k+1}=b'$, but not adjacent to $b_{k+2}$,
where $b_0b_1\cdots b_kb_{k+1}b'$ is a subpath of $H$.
In this case we replace the path $axyzb_0b_1,\cdots b_kb'$ with
$axzb_0b_1\cdots b_kyb'$ to obtain a new Hamilton cycle $H'$.
Note that none of $b=B_1,b_2,\ldots,b_k$ is a 2-vertex relative to $H'$ since none of these vertices is a
2-vertex with respect to $H$.
Moreover, whereas $H$ contained consecutive 2-vertices $x,y,z$, these have been replaced with $y$ and $b_k$ which are 2-vertices with respect to cycle $H'$.\par

Note that it is possible that $y$ is adjacent via an edge in $G_{int}$ to every vertex in $G$
(cf. Figure 2.4), but in that case
$\gamma (G)=1$.

\begin{figure}[!hbtp]
\begin{center}
\includegraphics[scale=.43]{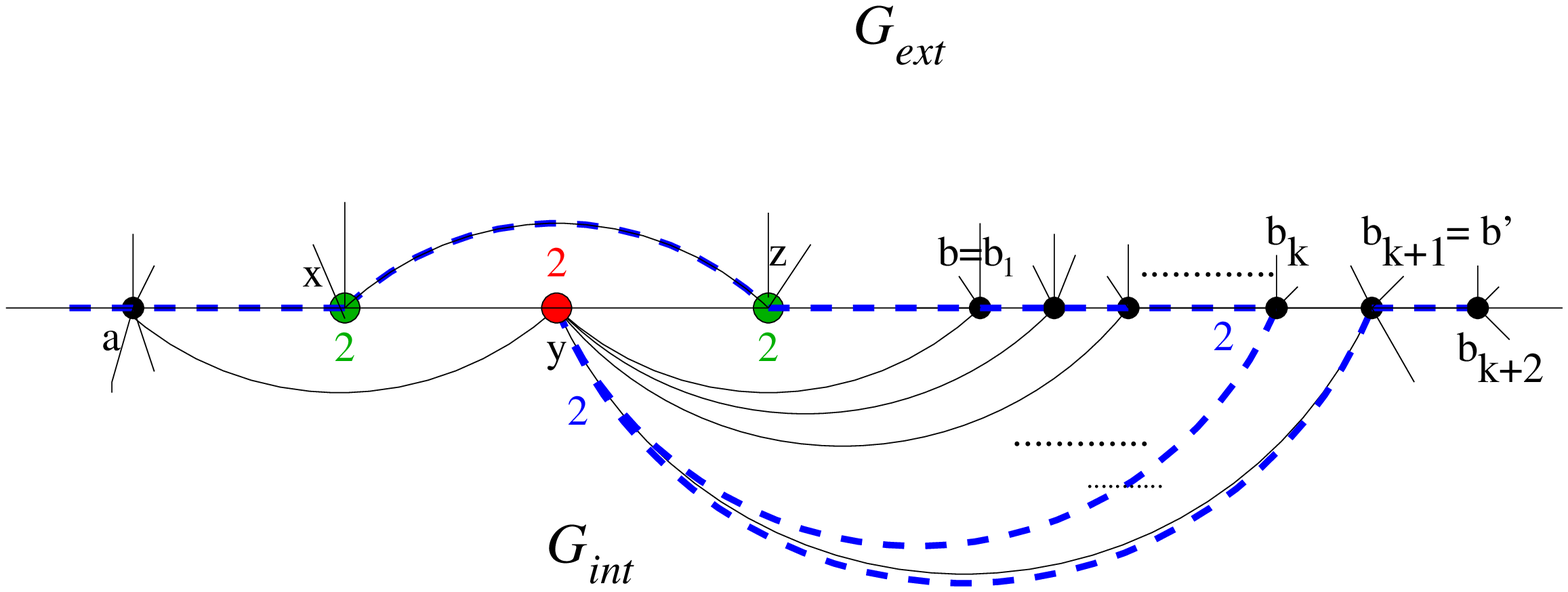} 
\centerline{\bf Figure 2.4}
\end{center}
\end{figure} 

\section{ $(H,A,B,O)$-graphs}\par

\medskip

In this section we introduce the concept of an $(H,A,B,O)$-graph and study some of its properties.\par

Let $G$ be a
(not necessarily planar)
graph on $n$ vertices with a Hamilton cycle $H = x_1x_2\cdots x_nx_1$.\par

A subgraph of $G$ consisting of a 3-cycle $x_ix_{i+1}x_{i+2}x_i$, where $\deg_G (x_{i+1})=2$, will be
denoted by
{\it $A$}.
A subgraph on four vertices consisting of the path $x_ix_{i+1}x_{i+2}x_{i+3}$, together with the edges
$x_ix_{i+2}$ and $x_{i+1}x_{i+3}$
such that $\deg_G(x_{i+1}) = \deg_{G}x_{i+2}=3$
will be denoted by {\it $B$}.
A third type of configuration, denoted by $O$, consists of a single edge with endvertices included
and such that if $x_i, i=1,2,$ is an endvertex of an $O$, then there is no edge of the form $x_{i-1}x_{i+1}$.
Further, if
an $A$ is immediately preceded and immediately succeeded by $O$s on the Hamilton cycle $H$, we will say that this $A$ is {\it isolated}.\par
\medskip
\noindent
{\bf Def.:}  An {\it $(A,B)$-string} is a connected maximal induced subgraph of $G$ consisting of $A$s and $B$s.
A {\it mixed} $(A,B)$-string is a string containing at least one $A$ and at least one $B$.\par
\medskip

%

We now suppose $G$, together with a Hamilton cycle $H$ contained in $G$, satisfies the following three assumptions.\par
\medskip
\noindent
(1) 
Suppose all edges in $E(G)-E(H)$ are of the form $x_ix_{i+2}$, ($\mod n$).
Such edges will be called {\it 2-chords}.\par
\medskip


\noindent
(2) Suppose further that $E(G)$ consists only of edges lying in the Hamilton cycle $H$ together with edges lying in either an
$A$ or a $B$.

\par
\medskip
\noindent
(3) Finally, assume that $G$ consists only of $A$s, $B$s and $O$s, where
any pair intersect in at most one vertex.\par
\medskip
Such a graph $G$, together with a Hamilton cycle $H$ in $G$, is called an
{\it (H,A,B,O)-graph}.
(The reader is referred to Figure 3.1 for an example of such a graph.)

\begin{figure}[!hbtp]\refstepcounter{figure}
\begin{center}
\includegraphics[scale=.45]{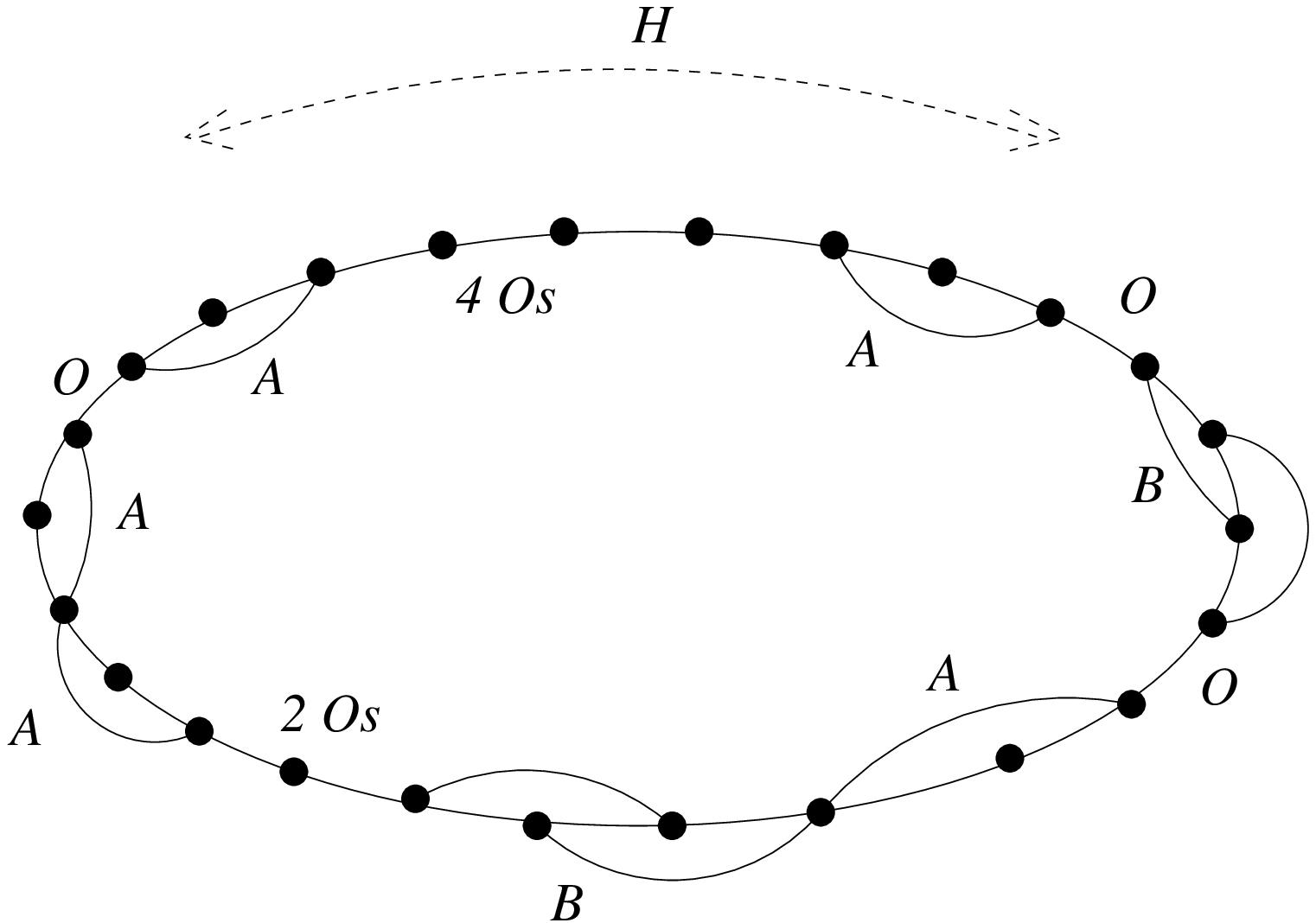} 
\centerline{\bf Figure 3.1}
\end{center}
\end{figure}

Note that an $(H,A,B,O)$-graph is planar.
%
%
Note also that it follows from Lemma 2.1 that every outerplanar triangulation with $\delta\ge 4$ which contains a Hamilton cycle must contain a spanning $(H,A,B,O)$-subgraph.\par

Let ${\cal K}$ denote the collection of all $(H,A,B,O)$-graphs with at least seven vertices and let
${\cal K}_{\lceil (n+1)/2\rceil}$ denote the subclass of ${\cal K}$ consisting of those members of ${\cal K}$ with
$n = |V(K)|$ vertices and in which the
Hamilton cycle $H$ possesses at least $(n+1)/2$ chords.
Let $K$ be a member of ${\cal K}_{\lceil (n+1)/2\rceil}$.
We now introduce seven operations,
called {\it reductions}.
In each of these reductions we replace a certain subgraph of $K\in {\cal K}_{\lceil (n+1)/2\rceil}$ with a certain smaller
subgraph so as
to produce a new graph $K'\in {\cal K}_{\lceil (n'+1)/2\rceil}$ for some $n'<n$.\par
We will also use a transformation which we will call a {\it switch} in which we replace a certain subgraph
of $K\in {\cal K}_{\lceil(n+1)/2\rceil}$ by another of the {\it same} length so as to produce a different graph which also
belongs to ${\cal K}_{\lceil(n+1)/2\rceil}$.
We then proceed to investigate how certain upper bounds on the domination number are affected when each of these
eight operations is carried out.
\par

In all switches and reductions to follow, the edges of $A$s and $B$s which are not edges of the Hamilton cycle $H$
may lie on either side of $H$, but of course no two such edges may cross since $K$ is planar.\par

\medskip
\noindent
{\bf Switch:}  Suppose $K\in {\cal K}_{\lceil (n+1)/2\rceil}$.
We modify $K$ by performing the replacment shown in
Figure 3.2, so as to form a new graph $K'$. \medskip

\begin{figure}[!hbtp] 
\begin{center}
\includegraphics[scale=.6]{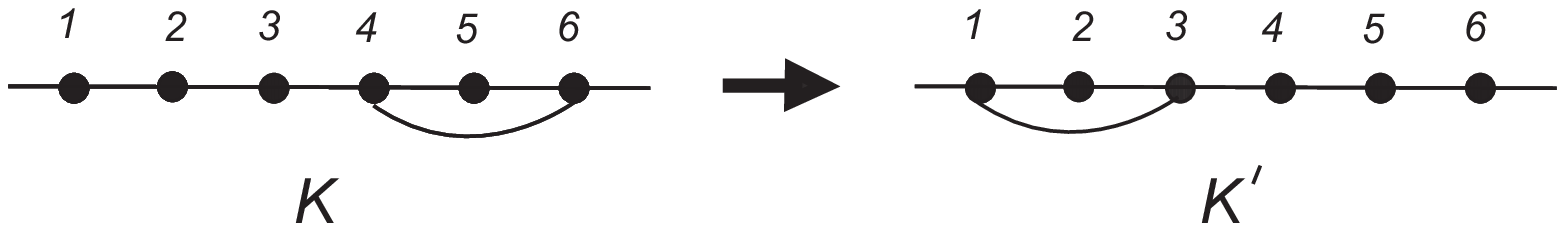} 
\centerline{\bf Figure 3.2}
\end{center}
\end{figure}

\noindent
{\bf Claim (a)}
$K'\in {\cal K}_{\lceil(n'+1)/2\rceil}$.
This is clear since $K$ and $K'$ have the same number of 2-chords.\par
\medskip
\noindent
{\bf Claim (b)}  $\gamma (K)\le \gamma (K')$.\par
\medskip
\noindent
{\it Proof of Claim (b):}
To prove Claim (b) we present a case-by-case analysis.
Suppose 1 and 6 are both in a minimum dominating set $D'$ for $K'$.
Then a third vertex is required to dominate 4 in $K'$. 
So then we may use 1,3 and 6 to dominate $K$.
Suppose 1 is in $D'$ and 6 is not in $D'$.
To dominate 5 in $K'$, one of 4 and 5 must be in $D'$.
Then use 1 and 4 to dominate $K$.
Suppose 1 is not in $D'$, but 6 is in $D'$.
Then one of 3 and 4 is in $D'$, so use 2 and 6 to dominate $K$.
Suppose neither 1 nor 6 is in $D'$.
Then two of 2, 3, 4 and 5 must be in $D'$,
so replace these two vertices with 2 and 4 to dominate $K$.
This proves Claim (b).\par
\medskip
We now turn our attention to reductions.
We formulate seven of these.
Again suppose $K\in{\cal K}_{\lceil (n+1)/2\rceil }$.
Let $t$ denote the number of 2-chords in $K$.

\medskip
\noindent
{\bf Reduction 1.}  Suppose $n\ge 11$.
Replace an $ABA$-segment in $K$ by a $B$-segment to obtain $K'$ as shown in Figure 3.3.

\medskip

\begin{figure}[!hbtp] 
\begin{center}
\includegraphics[scale=.6]{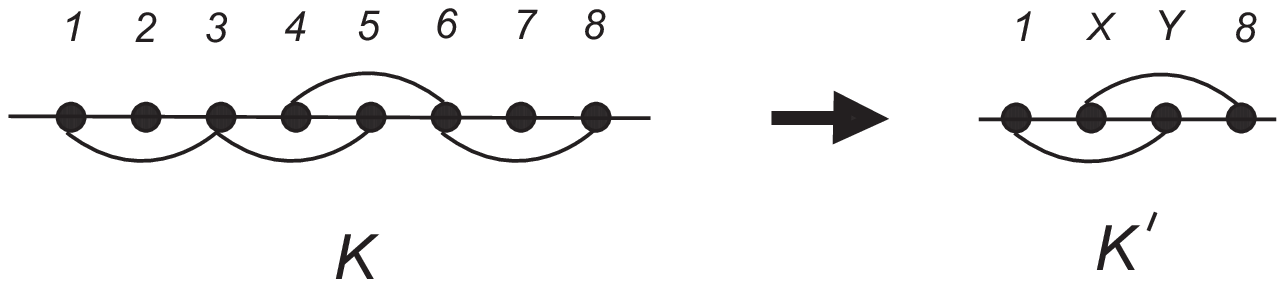} 
\centerline{\bf Figure 3.3}
\end{center}
\end{figure}

In this case, $n'=n-4$, $t'=t-2\ge (n+1)/2 -2=(n-3)/2=(n'+1)/2$,
and hence $K'\in {\cal K}_{\lceil(n'+1)/2\rceil}$.\par
\medskip
\noindent
{\bf Claim:}  
$\gamma (K)\le \gamma(K')+1$ and if $\gamma (K')\le |V(K')|/4$, then $\gamma (K)\le |V(K)|/4$.\par
\medskip
\noindent
{\it Proof of Claim:}
Let $D'$ be a minimum dominating set in $K'$.
Suppose 1 and 8 are in $D'$.
Then insert either 4 or 5 to obtain a dominating set $D$ for $K$ with $|D|=|D'|+1$.
Suppose 1 is in $D'$, but 8 is not in $D'$.
In this case we add 6 to $D'$ to dominate $K$.
Suppose 1 is not in $D'$, but 8 is in $D'$.
This case is symmetric to the preceding case.
Suppose neither 1 nor 8 is in $D'$.
Then one of $x$ and $y$ is in $D'$.
So we add 3 and 6 to dominate $K$.
In addition, if $\gamma (K')\le n'/4$, then $\gamma (K)\le n'/4 +1 = (n'+4)/4 = n/4$.
This proves the Claim.\par

\medskip
\noindent
{\bf Reduction 2.}  Suppose that $n\ge 11$.
Replace an $AA$-segment by a single vertex in $K$ to obtain $K'$ as shown in Figure 3.4.
\medskip

\begin{figure}[!hbtp] 
\begin{center}
\includegraphics[scale=.6]{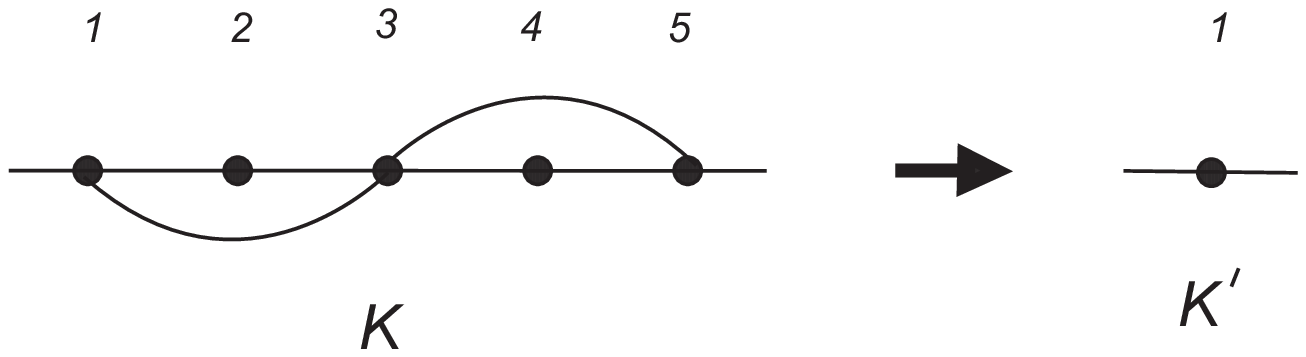} 
\centerline{\bf Figure 3.4}
\end{center}
\end{figure}

In this case, $n'=n-4$, $t'=t-2\ge (n+1)/2-2=(n-3)/2=(n'+1)/2$, and $K'\in {\cal K}_{\lceil (n'+1)/2\rceil}$.
\medskip
\noindent
{\bf Claim:}  $\gamma (K)\le \gamma (K')+1$ and if $\gamma (K')\le |V(K')|/4$, then
$\gamma (K)\le |V(K)|/4$.\par
\medskip
\noindent
{\it Proof of Claim:}  Let $D'$ be a minimum dominating set in $K'$.
If 1 is in $D'$, add vertex 5 to $D'$ to dominate $K$.
If 1 is not in $D'$, then
add 3 to $D'$ to dominate $K$.
Moreover, if $\gamma (K')\le n'/4$, then $\gamma (K)\le n'/4 + 1 = (n'+4)/4 = n/4$.
This proves the Claim.\par
\medskip\noindent
{\bf Reduction 3.}  Suppose $n\ge 11$.
Perform the transformation shown in Figure 3.5 to obtain $K'$.\par
\medskip

\begin{figure}[!hbtp] 
\begin{center}
\includegraphics[scale=.6]{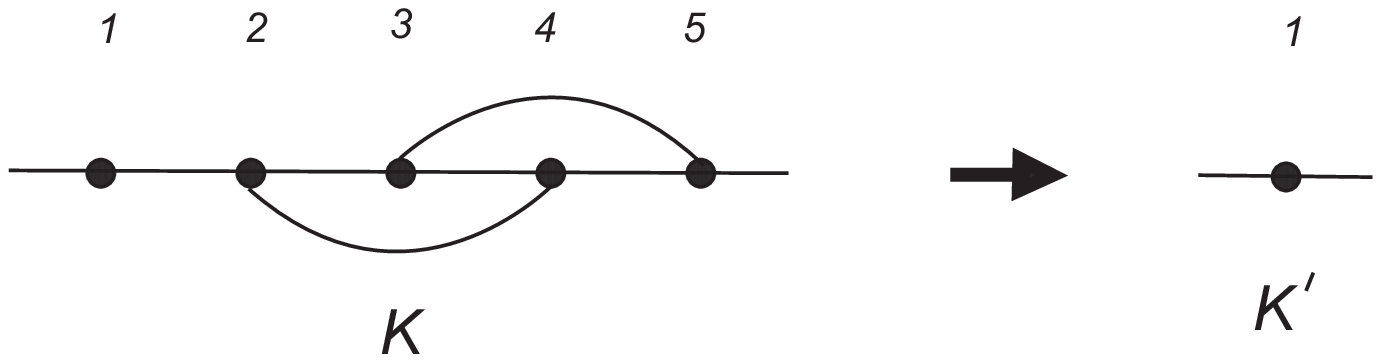} 
\centerline{\bf Figure 3.5}
\end{center}
\end{figure}

In this case, $n'=n-4$, $t'=t-2\ge (n+1)/2-2=(n-3)/2=(n'+1)/2$ and
$K'\in {\cal K}_{\lceil(n'+1)/2\rceil}$.\par
\medskip
\noindent
{\bf Claim:} 
$\gamma (K)\le \gamma (K')+1$.
Moreover, if $\gamma (K')\le |V(K')|/4$, then $\gamma (K)\le |V(K)|/4$.\par
\medskip
\noindent{\it Proof of Claim:}
Let $D'$ be a minimum dominating set in $K'$.
If 1 is in $D'$,
Then use 1 and 5 to dominate $K$.
If 1 is not in $D'$ and 1 is dominated from above, add 3 to dominate $K$.
If 1 is not in $D'$ and 1 is dominated from below,
add 2 to dominate $K$.
This completes the proof of the Claim.\par
\medskip
In the next two reductions, we reduce the total number of vertices by seven.
We want to show that if $\gamma (K')\le (2/7)|V(K')|$, then $\gamma (K)\le (2/7)|V(K)|$.\par
\medskip
\noindent
{\bf Reduction 4.}  Suppose $n\ge 14$.
Replace the 8-vertex configuration by a single vertex to obtain $K'$ as shown in Figure 3.6.\par
\medskip

\begin{figure}[!hbtp] 
\begin{center}
\includegraphics[scale=.6]{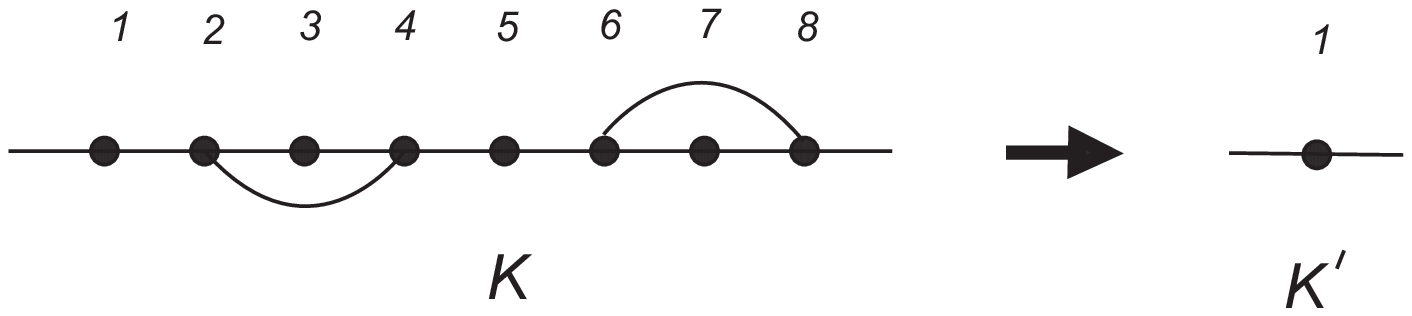} 
\centerline{\bf Figure 3.6}
\end{center}
\end{figure}

In this case $n'=n-7$, $t'=t-2\ge (n+1)/2-2=(n-3)/2>(n-6)/2=(n'+1)/2$ and $K'\in {\cal K}_{\lceil(n'+1)/2\rceil}$.\par
\medskip
\noindent
{\bf Claim:}
$\gamma (K)\le \gamma (K')+2$ and if $\gamma (K')\le (2/7)|V(K')|$, then
$\gamma (K)\le (2/7)|V(K)|$.
\par
\medskip
\noindent
{\it Proof of Claim:}
Let $D'$ be a minimum dominating set in $K'$.
If 1 is in $D'$, or 1 is not in $D'$, but is dominated from above, then add 4 and 8 to dominate $K$.
If 1 is not in $D'$, but is dominated from below in $K'$, then add 2 and 6 to dominate $K$.
This proves the Claim.\par
\par
\medskip
\noindent
{\bf Reduction 5.}  Suppose $n\ge 15$.
Replace the 13-vertex configuration by a 5-vertex configuration consisting of an $O$ and a $B$ to obtain $K'$ as shown in Figure 3.7.\par
\medskip

\begin{figure}[!hbtp] 
\begin{center}
\includegraphics[scale=.6]{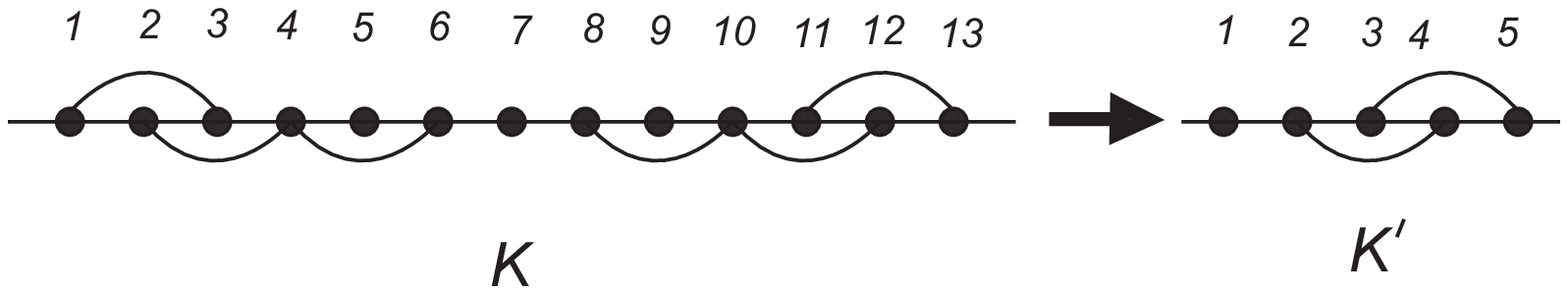} 
\centerline{\bf Figure 3.7}
\end{center}
\end{figure}

In this case, $n'=n-8\ge 7$ and $t'=t-4\ge (n+1)/2-4=(n'+1)/2$.
Hence $K'\in {\cal K}_{\lceil (n'+1)/2\rceil}$.\par
\medskip
\noindent
{\bf Claim:}
$\gamma (K)\le \gamma (K')+2$ and
$\gamma (K)\le 2|V(K)|/7$, if $\gamma(K')\le 2|V(K')|/7$.
\par
\medskip
\noindent{\it Proof of Claim:}
Let $D'$ be a dominating set in $K'$.
If 1 and $13\in D'$, then add 6 and 8 to dominate $K$.
If $1\in D'$, but $13\notin D'$,
then add 6 and 10 to dominate $K$.
\medskip
\noindent
If $1\notin D'$, but $13\in D'$, then one of 10, 11 or 12 belongs to $D'$ in order to dominate the vertex $10\in K'$.
Denote this vertex by $v$.
Then add 2, 4, and 8 to $D'-\{v\}$
to dominate $K$.
Finally suppose $1\notin D'$ and $13\notin D'$.
If both 1 and 13 are dominated by vertices from $\{10, 11, 12\}$, then $D'$ contains two vertices from $\{10,11,12\}$
which we will call $v$ and $w$.
Then add 3, 6, 8 and 11 to $D'-\{v,w\}$ to dominate $K$.
\par
So suppose at most one of 1 and 13 is dominated by vertices from $\{10,11,12\}$.
Then $D'$ contains at least one of 10, 11 and 12 which we will call $v$.
Then add 3, 6 and 10 to $D'-\{v\}$ to dominate $K$ if 13 is dominated by a vertex different from 10,11 and 12.
Otherwise, add 4, 8 and 11 to $D'-\{v\}$ to dominate $K$.
This completes the proof of the Claim.\par

\medskip
\noindent
{\bf Reduction 6.}  Suppose $n\ge 21$.
Replace the 15-vertex configuration by a single vertex to obtain $K'$ as shown in Figure 3.8.
Then again $K'\in {\cal K}_{\lceil(n'+1)/2\rceil}$.\par

In this case, $n'=n-14$, $t'=t-7\ge (n+1)/2-7=(n-13)/2=(n'+1)/2$.

\medskip
\noindent
{\bf Claim:}
$\gamma (K)\le \gamma (K')+4$ and
if $\gamma (K')\le (2/7)|V(K')|$, then $\gamma (K)\le (2/7)|V(K)|$.\par
\medskip

\begin{figure}[!hbtp] 
\begin{center}
\includegraphics[scale=.6]{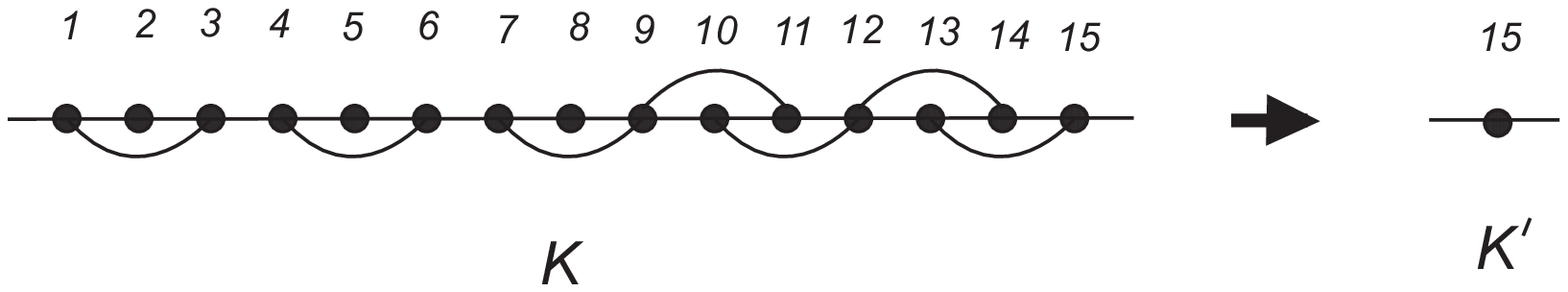} 
\centerline{\bf Figure 3.8}
\end{center}
\end{figure}

\noindent
{\it Proof of Claim:}
Let $D'$ be a minimum dominating set for $K'$.
If $15\notin D'$,
add 3, 6, 9 and 13 to $D'$.
This completes the proof of the Claim.
Our final reduction will be applied only in the final stages of the proof of Theorem 3.1.\par
\medskip
\noindent
{\bf Reduction 7.}  Suppose $n\ge 21.$
Replace the 12-vertex configuration by a single vertex to obtain $K'$ as shown in Figure 3.9.
\par
\bigskip

\begin{figure}[!hbtp] 
\begin{center}
\includegraphics[scale=.6]{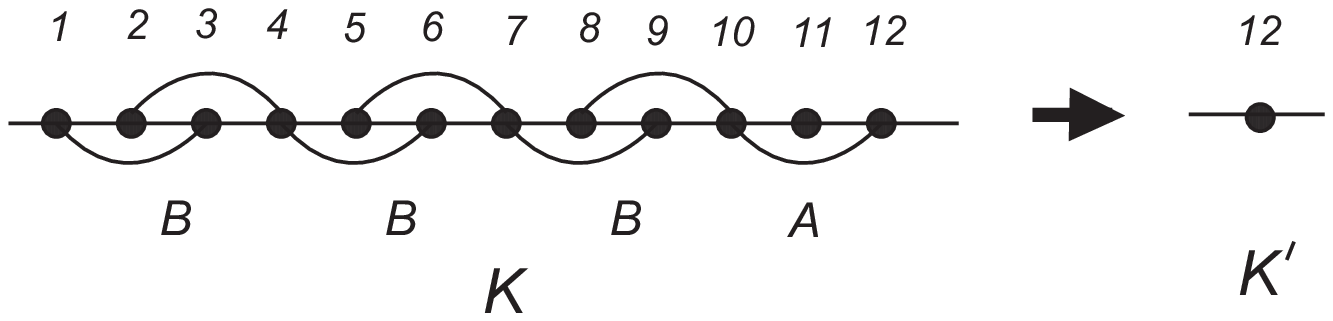} 
\centerline{\bf Figure 3.9}
\end{center}
\end{figure}
 
In this case, $|V(K')|=|V(K)|-11$.
\par
\medskip
\noindent
{\bf Claim:}
$\gamma (K)\le \gamma (K')+3$ and if $\gamma (K')\le
 2|V(K')|/7
$, then $\gamma (K)\le 
2|V(K)|/7
$.\par
\medskip
\noindent
{\it Proof of Claim:}
Let $D'$ be a minimum dominating set for $K'$.
If $12\in D'$, then add $1, 5, 7$ and $12$ to dominate $K$.
If $12\notin D'$, then add $2, 5$ and $10$ to dominate $K$.
Thus $\gamma (K)\le \gamma (K')+3\le
2|V(K')|/7
+3 =
2(|V(K)|-11)/7
+3\le
2|V(K)|/7
$
and the Claim is proved.
\par
\medskip
We are now prepared for our main result about $(H,A,B,O)$-graphs.\par
\medskip

\noindent
{\bf Theorem 3.1.}  If $K$ is an $(H,A,B,O)$-graph on $n$ vertices and $K$ has at least $(n+1)/2$ 2-chords,
then $\gamma (K)\le \lceil 2n/7\rceil$.\par
\medskip

\noindent
{\bf Proof:}  Let $K$ be an $(H,A,B,O)$-graph with 
$n$ vertices and at least $(n+1)/2$ 2-chords. Note that $K$ is simple.
It follows that $n\ge 5$ since 
$K$ has at least $(n+1)/2$ 2-chords. Suppose on the contrary that $K$ is a minimum counterexample.

\medskip

{\bf Claim 1}: {\sl $K$ has $n\ge 21$ vertices.}  \medskip

\noindent {\it Proof of Claim 1}: If not, suppose $5\le n\le 20$. Since $K$ has at least $ (n+1)/ 2 $ 2-chords, by the Pigenhole Principle, there are two
2-chords sharing a common endvertex $v$.
Then $v$ dominates five consecutive vertices on the 
Hamilton cycle.
But each vertex on the Hamilton cycle domintates at least three consecutive 
vertices.
Hence $\gamma(K)\le 1+\lceil (n-5)/3\rceil$. It can be easily checked that $1+\lceil (n-5)/3
\rceil\le \lceil 2n/7\rceil$ for $5\le n\le 20$.  So $\gamma(K)\le
\lceil 2n/7\rceil$, contradicting the assumption that $K$ is a counterexample. 
This completes the proof of Claim 1.
\medskip

In the following, we may therefore assume that $K$ has at least $n\ge 21$ vertices.
So Switch operations and Reductions 1-7  
can be applied to $K$.\medskip

{\bf Claim 2}: {\sl $K$ is composed of mixed $(A,B)$-strings only and each such mixed $(A,B)$-string $S$ satisfies (i), (ii) and (iii) below:

(i) $S$ starts and ends with an $A$; 

(ii) any two $A$s in $S$ are separated by at least two consecutive $B$s in $S$;

(iii) every such $S$
is separated from the next mixed $(A,B)$-string by exactly one $O$.}\medskip

\noindent{\it Proof of Claim 2}: We prove Claim 2 via five subclaims.
\medskip

(2.1) Every $(A,B)$-string in $K$ starts and ends with an $A$. \medskip

Suppose to the contrary that $S$ is an $(A,B)$-string which starts or ends with a $B$.
Applying Reduction 3, let 
$K'$ be the resulting $(H,A,B,O)$-graph which has four fewer vertices than $K$.
Since $K$ is a minimum counterexample, $\gamma(K')\le \lceil 2|V(K')|/7\rceil$.
By the claim following Reduction 3, we have $\gamma(K)\le \lceil 2n/7\rceil$,
which contradicts the assumption that $K$ is a minimum counterexample.
Subclaim (2.1) follows.\medskip

(2.2) No $(A,B)$-string of $K$ contains two consecutive $A$s. \medskip

Suppose to the contrary that $K$ has an $(A,B)$-string with two consecutive $A$s. 
Applying Reduction 2, let $K'$ be the new  
$(H,A,B,O)$-graph with $|V(K')|=n-4$ vertices. 
Then $\gamma(K')\le \lceil 2|V(K')|/7\rceil$ since $K$
is a minimum counterexample.
By the claim following Reduction 2, it follows that  $\gamma(K)\le \lceil 2n/7\rceil$, a contradiction
and subclaim (2.2) is proved. \medskip

(2.3) $K$ contains at least one mixed $(A,B)$-string and no mixed
$(A,B)$-string contains a segment of the form $ABA$.\medskip

Again by way of contradiction, suppose $K$ does not contain a mixed $(A,B)$-string.
Then by (2.1) and (2.2), all $(A,B)$-strings of
$K$ are isolated $A$s.
But then the number of 2-chords of $K$ is the number of 
$A$s in $K$, which is at most $n/3$, contradicting the assumption that $K$ has at least $(n+1)/2$
2-chords. Hence $K$ contains at least one mixed $(A,B)$-string.  

If $K$ has a mixed $(A,B)$-string with a segment of the form $ABA$, apply Reduction 1 to the segment
$ABA$.
Let $K'$ be the new $(H,A,B,O)$-graph.
Then $|V(K')|<|V(K)|$. Since $K$ is a minimum
counterexample, it follows that $\gamma(K')\le \lceil 2|V(K')|/7\rceil$.
By the claim following
Reduction 1, we have $\gamma(K)\le \lceil 2n/7\rceil$, contradicting the assumption that $K$ is a
counterexample.
Subclaim (2.3) thus follows. \medskip

(2.4)  For every vertex of degree 2 there is a 2-chord joining its neighbors; i.e., these three vertices form an $A$.
\medskip

By way of contradiction, suppose that $K$ has a vertex $v$ of degree 2 which is not contained in any $A$.
In other words, $v$ is the intersection of two consecutive $O$s. 
Let $w$ be a neighbor of $v$ in the Hamilton cycle.\par

First, assume that $w$ is also a vertex of degree 2.
Then $K$ contains three consective $O$s.
Using Switch operations, we may assume 
that one of $v$ and $w$, say $v$, has a neighbor in a mixed $(A,B)$-string $S$.
(We know such a mixed $(A,B)$-string exists by (2.3).)
If $v$ had no such neighbor, we could
apply Switch operations to bring the degree 2 vertices closer to the mixed $(A,B)$-string.
Hence, the three consecutive $O$s are followed by a mixed $(A,B)$-string $S$.
By (2.1) and (2.2), 
we know that $S=AB\cdots$; that is $S$ starts with an $A$ followed by a $B$. 
Apply a Switch operation to $OOOAB\cdots$ to obtain the segment $AOOOB\cdots$.
Then  we have a new $(H,A,B,O)$-graph $K'$ such that $|V(K')|=|V(K)|=n$ and $K'$ has a string
starting with a $B$.
By (2.1), $K'$ is not a counterexample and hence $\gamma(K')\le \lceil 2n/7\rceil$.
By the claim following the definition of a Switch operation, it follows that $\gamma(K)\le \gamma(K')\le 
\lceil 2n/7\rceil$ contradicting the assumption that $K$ is a minimum counterexample. 

So we may assume that both neighbors of $v$ have degree at least 3; i.e., $K$ contains two consecutive 
$O$s.
By (2.1), the two segments sharing vertices with the two consecutive $O$s
are $A$s.
If one of the $A$s is isolated, 
then apply Reduction 4 to obtain a new $(H,A,B,O)$-graph $K'$ such 
that $|V(K')|<|V(K)|$.
Since $K$ is a minimum counterexample, $\gamma(K')\le \lceil 2|V(K')|/7\rceil$.
By the Claim following Reduction 4,  we have $\gamma(K)\le \lceil 2n/7\rceil$,
contradicting the assumption that $K$ being a minimum counterexample.\par
So assume that neither of the $A$s is isolated.
Now apply Reduction 5.
As before, we then have $\gamma(K)\le \lceil 2n/7\rceil$, a contradiction, and (2.4) follows.\medskip
 
(2.5) Every $(A,B)$-string of $K$ is a mixed $(A,B)$-string.\medskip

Suppose to the contrary that $K$ contains an $(A,B)$-string which is not mixed.
Then by (2.1) and (2.2) this string must consist of an isolated $A$.
Among all such isolated $A$s, choose one that is closest to some mixed $(A,B)$-string.
It then follows by (2.3) and (2.4) that $K$ contains a segment of the form $AOAOABB$.
Apply Reduction 6 to $AOAOABB$ and let $K'$ be the new $(H,A,B,O)$-graph.
Then $|V(K')|<|V(K)|$.
Hence $\gamma(K')\le \lceil 2|V(K')|/7\rceil$ since $K$ is a minimum counterexample.
By the claim following Reduction 6, we have $\gamma(K)\le \lceil 2n/7\rceil$, a contradiction. \medskip

Combining (2.1)-(2.5), we can conclude that Claim 2 is proved.
 \medskip

{\bf Claim 3}: {\sl If $K$ contains at least six $B$s, then no three of these $B$s are consecutive.} \medskip

\noindent{\it Proof of Claim 3}: By way of contradiction, suppose that $K$ does contain three consecutive $B$s. 
Let $x_1$ be the number of $A$s and $x_2$ the number of $B$s.
Then $x_2\ge 6$ and since $K$ contains at least one mixed $(A,B)$-string, $x_1\ge 2$.\par

Let $y$ be the number of mixed $(A,B)$-strings. By Claim 2, it follows that $y \le (x_2-1)/2$ since at least one mixed $(A,B)$-string contains at least three consecutive $B$s.

Since $K$ contains three consecutive $B$s, $K$ contains a segment $ABBB$. 
Apply Reduction 7 to contract the segment $ABBB$ to a single vertex (deleting loops) and let $K'$ be the
resulting new $(H,A,B,O)$-graph. 
Then the number of 2-chords of $K'$ is $t'=(x_1-1)+2(x_2-3),$
and the number of vertices of $K'$ is $|V(K')| = 2(x_1-1) +3(x_2-3) +y,$
where $x_1-1$ and $x_2-3$ are the number of $A$s and $B$s in $K'$, respectively. 
But $x_2\ge 6$ and so it then follows that $x_2-3> (x_2-1)/2\ge y$.
Hence 

$${{t'}\over{|V(K')|}}={{(x_1-1)+2(x_2-3)}\over{2(x_1-1)+3(x_2-3)+y}} = 
   {{[(x_1-1)+(x_2-3)] +(x_2-3)}\over{2[(x_1-1)+(x_2-3)] + ((x_2-3)+y)}}
   > {{1}\over{ 2}}\ .$$
So the number of 2-chords in $K'$ is at least $(|V(K')|+1)/2$. Since $|V(K')|<|V(K)|$ and $K$ is a minimum counterexample,
$\gamma(K')\le \lceil 2|V(K')|/7 \rceil $.
By the claim following Reduction 7, we have $\gamma(K)
\le \lceil 2n/7\rceil$, contradicting the fact that $K$ is a counterexample.
This completes the proof of Claim 3.
\medskip 

By Claims 1-3, $S=ABBABBA\cdots ABBA$ for every mixed $(A,B)$-string $S$ in $K$ or else
$K$ has exactly five $B$s.

First assume that $S=ABBABBA\cdots ABBA$ for every mixed $(A,B)$-string in $K$.
Let $x$ be the total number of segments $ABB$ in $K$, and $y$ be the number of mixed $(A,B)$-strings. 
Then $n=8x+3y$.
Note that all vertices of the segment $ABB$
can be dominated by two vertices.
Since $x\ge y$, we have 
$$\gamma(K)\le 2x+y\le 2(7x+4y)/7 \le 2(8x+3y)/7=2n/7,$$
which contradicts the assumption that $K$ is a minimum counterexample.

On the other hand, if $K$ has exactly five $B$s,
then
$K=ABBAOABBBAO$ or else \hfill\break $K=ABBABBBAO$
by Claims 1 and 2.
Note that the domination number of $ABBA$ is 3 and
the domination number of $ABBBA$ is 3. Hence $\gamma(K)\le 6\le 2n/7$, a contradiction. 
This completes the
proof of the theorem. 

\section{Main Result}\par
\medskip
We shall have need of the following theorem due to Campos and Wakabayashi [CW].\par
\medskip
\noindent
{\bf Theorem 4.1:}  If $G$ is an maximal outerplanar graph with $n\ge 4$ vertices and with $t$ vertices of degree 2, then $\gamma (G)\le (n+t)/4$.\par
\medskip
We are now prepared for our main theorem.\par
\medskip
\noindent
{\bf Theorem 4.2:}  Let $G$ be a plane triangulation with $n$ vertices and $\delta (G)\ge 4$, and
suppose $G$ contains a Hamilton cycle $H$.
Then
$\gamma (G)\le \max\{\lceil 2n/7\rceil, \lfloor 5n/16\rfloor\}$.\par
\medskip
\noindent
{\bf Proof:}  If $\gamma (G)=1$, the result is trivial.
So suppose that $\gamma (G)\ge 2$. Since $\delta(G)\ge 4$, it follows that $n\ge 6$. 
First we use Lemma 2.1 to replace the given Hamilton cycle $H$ with one which has no three consecutive 2-vertices.\par

Suppose now that at least one of $G_{int}$ and $G_{ext}$ has no more than $n/4$ 2-chords,
say, without loss of generality, that $G_{int}$ has no more than $n/4$ 2-chords.
But then $G_{int}$ is a maximal outerplanar graph with $t\le n/4$ vertices of degree 2 and so by Theorem 4.1,
$\gamma (G)\le (n+t)/4\le (n+ n/4)/4= 5n/16$.\par

So suppose both $G_{int}$ and $G_{ext}$ have more than $n/4$ 2-chords.
Then the number of 2-chords of $G$ is larger than $n/4+n/4=n/2$.
In this case, we may apply Theorem 3.1 to conclude that $\gamma (G)\le \lceil (2n)/7\rceil$.
The proof of the theorem is complete. 
\par
\medskip

By a well-known theorem of Tutte [T], every 4-connected planar graph contains a Hamilton cycle, 
so the following result is an immediate corollary of Theorem 4.2.\medskip

\noindent
{\bf Corollary 4.3:}  Let $G$ be a 4-connected plane triangulation with $n$ vertices.
Then $\gamma (G)\le \max\{\lceil 2n/7\rceil, \lfloor 5n/16\rfloor\}$.\par


\section*{Concluding Remarks} 
Determination of the domination number seems to be difficult even for plane triangulations.
Matheson and Tarjan's conjecture has now been open for some eighteen years and there seems to be little significant progress toward
verifying their conjectured upper bound of $1/4|V(G)|$ for the general class of all planar triangulations.
The main result in the present paper makes some progress for those planar triangulations having minimum degree at least 4 and a Hamilton cycle.
(Note that this class is larger than the class of 4-connected planar triangulations.)
In particular, the reductions used in our main proof would seem to constitute a new approach at least for the graphs in this subclass.
Hopefully, they will prove to shed some light on the general case.

\end{document}